\documentclass[11pt]{amsart}
\usepackage[pdftex]{graphicx}
\usepackage{amsmath}
\usepackage{color}
\usepackage{amsthm}
\usepackage{multicol}

\usepackage{setspace} 

\usepackage{amssymb}

\usepackage{centernot}
\usepackage{bbm}

\usepackage{fancyvrb}
\usepackage{lipsum}

\theoremstyle{theorem}
\newtheorem{thm}{Theorem}[section]

\newtheorem{question}[thm]{Question}

\theoremstyle{definition}
\newtheorem{example}[thm]{Example}
\newtheorem{defn}[thm]{Definition}

\usepackage{mathtools}
\usepackage{enumitem}

\usepackage[margin=1.3in]{geometry}

\usepackage{verbatim}
\usepackage{url}
\usepackage{hyperref}
\usepackage{listings}
\usepackage{ragged2e}
\makeatletter

\renewcommand{\tocsection}[3]{%
  \indentlabel{\@ifnotempty{#2}{\bfseries\ignorespaces#1 #2\quad}}\bfseries#3}
\renewcommand{\tocsubsection}[3]{%
  \indentlabel{\@ifnotempty{#2}{\ignorespaces#1 #2\quad}}#3}

\newcommand\@dotsep{4.5}
\def\@tocline#1#2#3#4#5#6#7{\relax
  \ifnum #1>\c@tocdepth 
  \else
    \par \addpenalty\@secpenalty\addvspace{#2}%
    \begingroup \hyphenpenalty\@M
    \@ifempty{#4}{%
      \@tempdima\csname r@tocindent\number#1\endcsname\relax
    }{%
      \@tempdima#4\relax
    }%
    \parindent\z@ \leftskip#3\relax \advance\leftskip\@tempdima\relax
    \rightskip\@pnumwidth plus1em \parfillskip-\@pnumwidth
    #5\leavevmode\hskip-\@tempdima{#6}\nobreak
    \leaders\hbox{$\m@th\mkern \@dotsep mu\hbox{.}\mkern \@dotsep mu$}\hfill
    \nobreak
    \hbox to\@pnumwidth{\@tocpagenum{\ifnum#1=1\bfseries\fi#7}}\par
    \nobreak
    \endgroup
  \fi}
\AtBeginDocument{%
\expandafter\renewcommand\csname r@tocindent0\endcsname{0pt}
}
\def\l@subsection{\@tocline{2}{0pt}{2.5pc}{5pc}{}}
\makeatother

\author{Aaron Slobodin \\ \\Research Advisor: Dr. Sarah Mayes-Tang}

\address{Quest University Canada, 3200 University Blvd, Squamish, British Columbia, V8B 0N8}
\email{aaron.slobodin@questu.ca}

\begin{document}
\pagenumbering{gobble}
\title{Betti table stabilization of homogeneous monomial ideals}
\maketitle

\setcounter{section}{0}

\begin{abstract}
Given an homogeneous monomial ideal $I$, we provide a question- and example-based investigation of the stabilization patterns of the Betti tables shapes of $I^d$ as we vary $d$. We build off Whieldon's definition of the stabilization index of $I$, Stab$(I)$, to define the stabilization sequence of $I$, StabSeq$(I)$, and use it to explore changes in the shapes of the Betti tables of $I^d$ as we vary $d$. We also present the stabilization indices and sequences of the collection of ideals $\{I_{n}\}$ where $I_{n}=(a^{2n}b^{2n}c^{2n},b^{4n}c^{2n},a^{3n}c^{3n},a^{6n-1}b)\subseteq \mathbbm{k}[a,b,c]$.
\end{abstract}

\section{Introduction}
There is little known about the relationship between the Betti tables $\beta(I^d)$ of $I^d$ as we vary $d$. Elena Guardo and Adam Van Tuyl compute the Betti numbers of homogeneous complete intersection ideals, in \cite{EG}. In \cite{GW}, Whieldon proved that the shapes of the Betti tables of equigenerated ideals will stabilize, and provides a conjecture for the formula for the power of stabilization for edge ideals $I_{G}$. By ``shape" of the Betti tables, of a given ideal $I$, we mean the arrangement of their nonzero entries and by ``stabilize" we mean that there exists some $D$ for which the Betti tables of $I^d$ will share the same shape as $I^D$ for all $d \geq D$. Though Whieldon's proof ensures the $D$ for which we see a stabilized Betti table shape of $\{I^d\}$ is finite, there is currently no formula to predict the value of $D$.

The lack of a generalized formula for the power of $I$, $D$, for which we see a stabilized shape of $\beta (I^d)$ presents a problem for computational investigations. Though we might be confident, having seen no evidence to the contrary, that the Betti tables of $\{I^d\}$ have stabilized in shape at $D$, we often have no way of knowing. Though this lack of knowledge is problematic, it also opens the door to exciting questions. Such questions include: If the Betti tables of $I^D$ through $I^{D+x}$ share the same shape, for some $x$, can we feel confident that the shape of the Betti tables of $\{I^d\}$ have stabilized by $D$? If so, what are the bounds on this $x$? In what ways can we characterize the Betti table shape stabilization patterns of a given ideal? Are there ideals for which we do not see a stabilized Betti table shape until $I^d$ where $d$ is arbitrarily large?

In this paper we focus on the Betti table shape stabilization patterns of homogeneous monomial ideals and provide an example-driven narrative to unearth the potential of these questions. 

\hyperref[B:Background]{Section~\ref*{B:Background}} provides background on Betti numbers, Betti tables, and the stabilization index. In this section we also define the \textit{stabilization sequence} of a given ideal, and present Whieldon's result guaranteeing Betti table shape stabilization. 

\hyperref[Q: Questions and Results]{Section~\ref*{Q: Questions and Results}} is framed around questions and examples concerning patterns in the stabilization indices and sequences of homogeneous monomial ideals. In this section we also define a \textit{linearly connected family} of ideals and investigate the linearly connected family of ideals $\{I_{n}\}$ where $I_{n}=(a^{2n}b^{2n}c^{2n},b^{4n}c^{2n},a^{3n}c^{3n},a^{6n-1}b)\subseteq \mathbbm{k}[a,b,c]$ which demonstrates nice relationships between the stabilization indices and sequences of $I_{n}$ as we vary $n$. 

\hyperref[Appendix A]{Section~\ref*{Appendix A}} presents the Macaulay2 \cite{M2} commands used in this paper. In this section we also provide the Macaulay2 \cite{M2} session, with commentary, used in Example\autoref{E: Example 2.1} and Example\autoref{E:Example 2.3}.

\hyperref[Appendix B]{Section~\ref*{Appendix B}} provides the source code of the Macaulay2 \cite{M2} package \texttt{StabSeq.m2}. This package, written by the author, is used to produce the stabilization sequences prestented in this paper.

\section{Background} \label{B:Background}
Let $R=\mathbbm{k}[x_{1},x_{2},...,x_{n}]$, where $\mathbbm{k}$ is an algebraically closed field of characteristic zero. For $j \in Z$ denote $R(-j)$ to be a graded R-module such that $R(-j)_{i}=R_{i-j}$. We say that $R(-j)$ is shifted $j$ degrees. An element of degree $i$ in $R$ would have degree $i+j$ in $R(-j)$. Following the convention in \cite{Peeva2011} the element $1 \in R(-j)$ has degree $j$ and is called the \textit{1-generator} of $R(-j)$. 

\subsection{Minimal graded free resolutions} \label{B: Minimal graded free resolutions}
Let $I \subseteq R$ be a homogeneous monomial ideal. Given $I$ we can construct a \textit{free resolution of $I$} 
$$
0 \xrightarrow{d_{\ell+1}} F_{\ell} \xrightarrow{d_{\ell}} F_{\ell-1} \xrightarrow{d_{\ell-1}} \cdots \xrightarrow{d_{2}} F_{1} \xrightarrow{d_{1}} F_{0} \xrightarrow{d_{0}} I \xrightarrow{d} 0.
$$
such that: \begin{enumerate}
\item for all $i$, $F_{i}$ is a free finitely generated graded R-module of the form
$$
F_{i} = \oplus_{j \in \mathbbm{Z}}R(-j)^{\beta_{i,j}(I)}
$$
for some graded Betti numbers $\beta_{i,j}(I)$,
\item the sequence of homomorphism is \textit{exact}, that is Ker$(d_{i})=\text{Im}(d_{i+1})$, for all $i$, and Ker$(d)=\text{Im}(d_{0})$. 
\end{enumerate}

\begin{defn}
A free resolution is \textit{minimal} if $d_{i+1}(F_{i+1})\subseteq(x_{1},x_{2},...,x_{n})F_{i}$, for all $i$. 
\end{defn} 

Informally, this condition means that the differential matrix, $D_{i}$, given by the map $d_{i}$, contains no nonzero constants. 

\begin{example}
The homomorphism
$$
R(-5) \oplus R(-4) \xrightarrow{
\begin{pmatrix}
1 & 0\\
0 & x^2\\
\end{pmatrix}} R(-5) \oplus R(-2)
$$
is not minimal since the differential matrix of the map from $R(-5) \oplus R(-4)$ to $R(-5) \oplus R(-2)$ contains the nonzero constant 1. Note that the map from $R(-5)$ to $R(-5)$ is trivial, the map being multiplication by the identity. In general, we can construct a minimal free resolution from a free resolution by removing all trivial homomorphism. 
\end{example}

\begin{defn}
A free resolution is \textit{graded} if:
\begin{enumerate}
\item each R-module $F_{i}$ is graded,
\item $R_{j}F_{i} \subseteq F_{i+j}$ for all $i, j \in \mathbbm{Z}$,
\item the map $d_{i}$, for all $i$, is a homomorphism of degree 0.
\end{enumerate}
\end{defn}
We say that a homomorphism $d_{i}:F_{i} \rightarrow F_{i-1}$ has \textit{degree} $j$ if deg($d_{i}(n))=j+\text{deg}(n)$ for each $n \in F_{i}$.

\begin{example}
Consider the homomorphism
$$
R(-3) \oplus R(-2) \xrightarrow{
\begin{pmatrix}
x_{1}^3 & x_{2}^2\\ 
\end{pmatrix}} R.
$$
Let $(f,g) \in R(-3) \oplus R(-2)$. Under this homomorphism $(f,g) \mapsto fx_{1}^3+gx_{2}^2$. Let $f$ and $g$ have degrees $a_{1}$ and $a_{2}$, respectively. The deg$(fx_{1}^3)=a_{1}+3$ and deg$(gx_{2}^2)=a_{2}+2$ in $R$. Since the degrees of $f \in R(-3)$ and $g \in R(-2)$ are shifted by 3 and 2, respectively, the deg$(f)=a_{1}+3$ in $R(-3)$ and deg$(g)=a_{2}+2$ in $R(-2)$. Since deg$(fx_{1}^3)=\text{deg}(f)$ and deg$(gx_{2}^2)=\text{deg}(g)$ the homomorphism has degree 0.

\end{example}

We now provide a step by step construction of a minimal graded free resolution in the following example.

\begin{example} \label{E: Example 2.1}
Let $I=(x_{1}x_{2}^2,x_{1}x_{3}^2,x_{2}^3,x_{1}^3) \subseteq \mathbbm{k}[x_{1},x_{2},x_{3}]$. We begin our minimal graded free resolution of $I$ with
$$
F_{0} \xrightarrow{d_{0}} I \xrightarrow{d} 0.
$$

The map $d$ is surjective, therefore the generators of $I$, all of which have degree 3, generate Ker$(d)$. Since Ker$(d)=\text{Im}(d_{0})$, the generators of $I$ also generate Im$(d_{0})$. We therefore set $F_{0}=R(-3) \oplus R(-3) \oplus R(-3) \oplus R(-3)$. Denote $f_{1},f_{2},f_{3},f_{4}$ to be the homogeneous 1-generator of the 4 $R$ modules of $F_{0}$. Note that, for all $i$, deg$(f_{i})=3$. We now define $d_{0}$ by 
$$
\begin{tabular}{cccc}
$f_{1} \mapsto x_{1}x_{2}^2,$ & $f_{2} \mapsto x_{1}x_{3}^2,$ & $f_{3} \mapsto x_{2}^3,$ & $f_{4} \mapsto x_{1}^3$.\\
\end{tabular}
$$
Our resolution therefore begins
$$
F_{1} \xrightarrow{d_{1}} R^4(-3) \xrightarrow{
\begin{pmatrix}
x_{1}x_{2}^2 & x_{1}x_{3}^2 & x_{2}^3 & x_{1}^3\\
\end{pmatrix}
} I \rightarrow 0.
$$

Since Ker$(d_{0})=\text{Im}(d_{1})$, we can determine the generating set of Im$(d_{1})$ by determining the generating set of Ker$(d_{0})$. Let $\alpha_{1}f_{1}+\alpha_{2}f_{2}+\alpha_{3}f_{3}+\alpha_{4}f_{4} \in \text{Ker}(d_{0})$ where, for all $i$, $\alpha_{i} \in R$. To determine the generators of Ker$(d_{0})$, we must solve the equation
$$
\alpha_{1}x_{1}x_{2}^2+\alpha_{2}x_{1}x_{3}^2+\alpha_{3}x_{2}^3+\alpha_{4}x_{1}^3=0
$$
for $(\alpha_{1},\alpha_{2},\alpha_{3},\alpha_{4})$. The generating solutions to the above equation are
$$
\begin{tabular}{cccc}
$\sigma_{1}=(-x_{2},0,x_{1},0)$, & $\sigma_{2}=(x_{1}^2,0,0,-x_{2}^2)$, & $\sigma_{3}=(0,x_{1}^2,0,-x_{3}^2)$, & $\sigma_{4}=(-x_{3}^2,x_{2}^2,0,0)$ \\
\end{tabular}
$$

Thus
$
\begin{tabular}{cccc}
$-x_{2}f_{1}+x_{1}f_{3},$ & $x_{1}^2f_{1}-x_{2}^2f_{4},$ & $x_{1}^2f_{2}-x_{3}^2f_{4},$ & $-x_{3}^2f_{1}+x_{2}^2f_{2}$ \\
\end{tabular}
$
generate Ker$(d_{0})$. Their degrees are $4,5,5,5$ respectively. We therefore set $F_{1}=R(-4) \oplus R(-5) \oplus R(-5) \oplus R(-5)$. Denote $g_{1}, g_{2}, g_{3}, g_{4}$ to be the homogeneous 1-generator of the 4 R-modules of $F_{1}$, respectively. Note that deg$(g_{1})=4$ and deg$(g_{i})=5$ for $2 \leq i \leq 4$, . We now define $d_{1}$ by 
$$
\begin{tabular}{cccc}
$g_{1} \mapsto -x_{2}f_{1}+x_{1}f_{3},$ & $g_{2} \mapsto x_{1}^2f_{1}-x_{2}^2f_{4},$ & $g_{3} \mapsto x_{1}^2f_{2}-x_{3}^2f_{4},$ & $g_{4} \mapsto -x_{3}^2f_{1}+x_{2}^2f_{2}$\\
\end{tabular}
$$
with the resulting next step of our resolutions
$$
F_{2} \xrightarrow{d_{2}}
R(-4) \oplus R^3(-5) \xrightarrow{
\begin{pmatrix}
-x_{2} & x_{1}^2 & 0 & -x_{3}^2\\
0 & 0 & x_{1}^2 & x_{2}^2\\
x_{1} & 0 & 0 & 0\\
0 & -x_{2}^2 & -x_{3}^2 & 0\\
\end{pmatrix}} R^4(-3) \xrightarrow{
\begin{pmatrix}
x_{1}x_{2}^2 & x_{1}x_{3}^2 & x_{2}^3 & x_{1}^3\\
\end{pmatrix}
} I \rightarrow 0.
$$

Recall Ker$(d_{1})=\text{Im}(d_{2})$. Let $\beta_{1}g_{1}+\beta_{2}g_{2}+\beta_{3}g_{3}+\beta_{4}g_{4} \in \text{Ker}(d_{1})$ where, for all $i$, $\beta_{i} \in R$. To determine the generators of Ker$(d_{1})$, we must solve the equation
$$
\beta_{1}(-x_{2}f_{1}+x_{1}f_{3})+\beta_{2}(x_{1}^2f_{1}-x_{2}^2f_{4})+\beta_{3}(x_{1}^2f_{1}-x_{2}^2f_{4})+\beta_{4}(-x_{3}^2f_{1}+x_{2}^2f_{2})=0
$$
for $(\beta_{1},\beta_{2},\beta_{3},\beta_{4})$. The generating solution to the above equation is
$$
\eta_{1}=(0,-x_{3}^2,-x_{2}^2,x_{1}^2).
$$

Thus $-x_{3}^2g_{2}-x_{2}^2g_{3}+x_{1}^2g_{4}$ generates Ker$(d_{1})$. The degree of the generator of Ker$(d_{1})$ is 7. We therefore set $F_{2}=R(-7)$. Denote $h_{1}$ to be the homogeneous 1-generator of $R(-7)$. Note that deg($h_{1})=3$. We now define $d_{2}$ by
$$
h_{1} \mapsto -x_{3}^2g_{2}-x_{2}^2g_{3}+x_{1}^2g_{4}
$$
with the resulting next step of our resolutions
$$
F_{3} \xrightarrow{d_{3}}
R(-7) \xrightarrow{
\begin{pmatrix}
0\\
-x_{3}^2\\
-x_{2}^2\\
x_{1}^2\\
\end{pmatrix}
}
R(-4) \oplus R^3(-5) \xrightarrow{
\begin{pmatrix}
-x_{2} & x_{1}^2 & 0 & -x_{3}^2\\
0 & 0 & x_{1}^2 & x_{2}^2\\
x_{1} & 0 & 0 & 0\\
0 & -x_{2}^2 & -x_{3}^2 & 0\\
\end{pmatrix}}
$$
$$
R^4(-3) \xrightarrow{
\begin{pmatrix}
x_{1}x_{2}^2 & x_{1}x_{3}^2 & x_{2}^3 & x_{1}^3\\
\end{pmatrix}
} I \rightarrow 0.
$$

Recall Ker$(d_{2})=\text{Im}(d_{3})$. Let $\gamma_{1}h_{1} \in \text{Ker}(d_{2})$, where $\gamma_{1} \in R$. To determine the generators of Ker$(d_{2})$, we must solve the equation
$$
\gamma_{1}h_{1}=0
$$

The generating solution to the above equation is $\gamma_{1}=0$. Thus 0 generates Ker$(d_{2})$ and Im$(d_{3})$. We conclude by setting $F_{3}=0$. We now obtain our complete minimal graded free resolution of $I$.

$$
0 \rightarrow
R(-7) \xrightarrow{
\begin{pmatrix}
0\\
-x_{3}^2\\
-x_{2}^2\\
x_{1}^2\\
\end{pmatrix}
}
R(-4) \oplus R^3(-5) \xrightarrow{
\begin{pmatrix}
-x_{2} & x_{1}^2 & 0 & -x_{3}^2\\
0 & 0 & x_{1}^2 & x_{2}^2\\
x_{1} & 0 & 0 & 0\\
0 & -x_{2}^2 & -x_{3}^2 & 0\\
\end{pmatrix}}
$$
$$
R^4(-3) \xrightarrow{
\begin{pmatrix}
x_{1}x_{2}^2 & x_{1}x_{3}^2 & x_{2}^3 & x_{1}^3\\
\end{pmatrix}
} I \rightarrow 0.
$$
\end{example}

By construction our resolution is graded and free. Our graded free resolution is also minimal, given that no nonzero constants are contained in any of the differential matrices $D_{i}$. All remaining minimal graded free resolutions in this paper were computed in Macaulay2 \cite{M2}. See \hyperref[Appendix]{Section~\ref*{Appendix A}} for Macaulay2 \cite{M2} commands used in this paper. 

\subsection{Graded Betti numbers and Betti numbers} \label{Graded Betti numbers and Betti numbers}

Below we use the same notation as in \hyperref[B: Minimal graded free resolutions]{Section~\ref*{B: Minimal graded free resolutions}}.

\begin{defn}
The \textit{graded Betti numbers} of $I$ are the numbers $$\beta_{i,j}(I)= \text{number of copies of}\ R(-j)\ \text{in}\ F_{i}.$$
\end{defn}
\begin{defn}
The $i^{th}$ \textit{total Betti number} of $I$ is the sum of the graded Betti numbers in homological degree $i$: $$\beta_{i}(I)= \sum_{j \in \mathbbm{Z}}\beta_{i,j}(I).$$
\end{defn}
\begin{example} \label{E:Example 2.3}
In Example\autoref{E: Example 2.1}, we saw that the minimal graded free resolution of $I=(x_{1}x_{2}^2,x_{1}x_{3}^2,x_{2}^3,x_{1}^3) \subseteq \mathbbm{k}[x_{1},x_{2},x_{3}]$ was
$$
0 \rightarrow R(-7)  \rightarrow R(-4) \oplus R^3(-5) \rightarrow R^4(-3) \rightarrow I \rightarrow 0.
$$

In the above minimal graded free resolution of $I$ we have that $F_{0}=R^4(-3)$, $F_{1}=R(-4) \oplus R^3(-5)$, and $F_{2}=R(-7)$. The graded Betti numbers of $I$ are therefore

$$
\begin{tabular}{cccc}
$\beta_{0,3}(I)=4$ & $\beta_{1,4}(I)=1$ & $\beta_{1,5}(I)=3$ & $\beta_{2,7}(I)=1$\\
\end{tabular}
$$
and the total Betti numbers of $I$ are
$$
\begin{tabular}{ccc}
$\beta_{0}(I)=4$ & $\beta_{1}(I)=4$ & $\beta_{2}(I)=1$.\\
\end{tabular}
$$
\end{example}

\subsection{Betti tables}
\begin{defn}
Once we have determined the graded Betti numbers of a given ideal, we can express them in a \textit{Betti table}. The Betti tables in this paper are styled as computed in Macaulay2 \cite{M2}.  The $i^{th}$ column and $j^{th}$ row in a Betti table corresponds to the graded Betti number $\beta_{i,i+j}$. A graded Betti number equal to 0 in a Betti table will be denoted by a $\cdot$ .
\end{defn}
\begin{example} \label{E:Example 2.5}
The Betti table of the ideal in Example\autoref{E:Example 2.3} would be expressed as follows.
$$
\begin{tabular}{ccccc}
$I$ & & & $I$ &\\
 & 
 \begin{tabular}{cccc}
- & 0 & 1 & 2\\
total: & $\beta_{0}(I)$ & $\beta_{1}(I)$ & $\beta_{2}(I)$\\
1: & $\beta_{0,1}(I)$ & $\beta_{1,2}(I)$ & $\beta_{2,3}(I)$\\
2: & $\beta_{0,2}(I)$ & $\beta_{1,3}(I)$ & $\beta_{2,4}(I)$\\
3: & $\beta_{0,3}(I)$ & $\beta_{1,4}(I)$ & $\beta_{2,5}(I)$\\
4: & $\beta_{0,4}(I)$ & $\beta_{1,5}(I)$ & $\beta_{2,6}(I)$\\
5: & $\beta_{0,5}(I)$ & $\beta_{1,6}(I)$ & $\beta_{2,7}(I)$\\
\end{tabular}
& $\rightarrow$ & &
 \begin{tabular}{cccc}
- & 0 & 1 & 2\\
total: & 4 & 4 & 1\\
1: & $\cdot$ & $\cdot$ & $\cdot$\\
2: & $\cdot$ & $\cdot$ & $\cdot$\\
3: & 4 & 1 & $\cdot$\\
4: & $\cdot$ & 3 & $\cdot$\\
5: & $\cdot$ & $\cdot$ & 1\\
\end{tabular}\\
\end{tabular}
$$
\end{example}
The Betti tables in the rest of this paper will be shifted to only show nonzero Betti numbers in the resolution of $I^d$. The Betti table in Example\autoref{E:Example 2.5} would therefore be expressed as follows.
$$
\begin{tabular}{cc}
$I$ & \\
&
 \begin{tabular}{cccc}
- & 0 & 1 & 2\\
total: & 4 & 4 & 1\\
3: & 4 & 1 & $\cdot$\\
4: & $\cdot$ & 3 & $\cdot$\\
5: & $\cdot$ & $\cdot$ & 1\\
\end{tabular}\\
\end{tabular}
$$
\subsection{Betti table shape}
\begin{defn}
We say that the Betti tables of $I^x$ and $I^y$ share the same \textit{shape} if there exists an integer $r$ such that, for all $i$ and $j$,
$$
\beta_{i,j+rx}(I^x) \neq 0 \Longleftrightarrow \beta_{i,j+ry}(I^y) \neq 0.
$$ 
\end{defn}

Informally, this condition means that the Betti tables of $I^x$ and $I^y$ exhibit the same structure of nonzero entries. 

Notice that the $r$ in this definition corresponds to the degree of the lowest degree generator of $I$: if $g$ is the generator of lowest degree $r$ in $I$, $g^x$ and $g^y$ will be the generators of lowest degree in $I^x$ and $I^y$, respectively, and correspond to the least nonzero entry in the $0^{th}$ column of the Betti tables $I^x$ and $I^y$, respectively. It follows that the least nonzero entry in the $0^{th}$ column of the Betti table of $I^x$ and $I^y$ will therefore occur, respectively, in row $rx$ and $ry$. Therefore $r$, in the condition
$$
\beta_{i,j+rx}(I^x) \neq 0 \Longleftrightarrow \beta_{i,j+ry}(I^y) \neq 0
$$ 
can be considered as a shift of the rows of the Betti tables of $I^x$ and $I^y$, dependant upon the generator of lowest degree of $I$. 

The following example will provide clarity to this concept of Betti table shape.

\begin{example}  \label{E:Example 2.7}
Let $I=(x_{1}^2,x_{2}^2,x_{3}^2,x_{4}^2) \subseteq \mathbbm{k}[x_{1},x_{2},x_{3},x_{4}]$. The Betti tables $\beta (I)$ and $\beta (I^2)$ are given by:
$$
\begin{tabular}{ccccc}
$I$ & & & $I^2$ & \\
 &
\begin{tabular}{ccccc}
- & 0 & 1 & 2 & 3\\
total: & 4 & 6 & 4 & 1\\
2: & 4 & $\cdot$ & $\cdot$ & $\cdot$\\
3: & $\cdot$ & 6 & $\cdot$ & $\cdot$\\
4: & $\cdot$ & $\cdot$ & 4 & $\cdot$\\
5: & $\cdot$ & $\cdot$ & $\cdot$ & 1\\
\end{tabular}
& & &
\begin{tabular}{ccccc}
- & 0 & 1 & 2 & 3\\
total: & 10 & 20 & 15 & 4\\
4: & 10 & $\cdot$ & $\cdot$ & $\cdot$\\
5: & $\cdot$ & 20 & $\cdot$ & $\cdot$\\
6: & $\cdot$ & $\cdot$ & 15 & $\cdot$\\
7: & $\cdot$ & $\cdot$ & $\cdot$ & 4\\
\end{tabular}\\
\end{tabular}
$$

We say that the Betti tables of $I^1$ and $I^2$ share the same shape since they both exhibit the same diagonal structure of nonzero entries. In the above Betti tables we see that 
$$
\begin{tabular}{cccc}
$\beta_{0,2}(I) \Longleftrightarrow \beta_{0,4}(I^2)$, & $\beta_{1,4}(I) \Longleftrightarrow \beta_{1,6}(I^2)$, & $\beta_{2,6}(I) \Longleftrightarrow \beta_{2,8}(I^2)$, & $\beta_{3,8}(I) \Longleftrightarrow \beta_{3,10}(I^2),$\\
\end{tabular}
$$
Thus, for all $i$ and $j$, 
$$\beta_{i,j+2(1)}(I^1) \neq 0 \Longleftrightarrow \beta_{i,j+2(2)}(I^2) \neq 0.$$
\end{example}
Though the Betti tables in Example\autoref{E:Example 2.7} shared the same shape, not all Betti tables of different powers of a given ideal do. There are many simple examples for which this is true, several of which will be explored in the rest of this paper. We will now develop a language to better handle different Betti table shapes.

\subsection{Betti table shape stabilization}
\begin{defn} \textnormal{(Definition 1.2 in \cite{GW})}
We say that an ideal $I=(f_{0},f_{1},...,f_{k}) \subseteq R$ is \textit{equigenerated in degree r} if $deg(f_{i})=r$ for all $f_{i}$.
\end{defn}
The following result by Whieldon shows that, given $I \subseteq R$, where $I$ is an equigenerated ideal of degree $r$, the Betti tables of $I^d$ will \textit{stabilize}, in shape, for all $d$ greater than some index $D$. That is, though the Betti tables of $I^d$, where $d < D$ may differ in shape, all Betti tables of $I^d$ where $d \geq D$ will share the same shape as $I^D$. 

\begin{thm} \label{T: Theorem 2.9}
\textnormal{(Theorem 4.1 in \cite{GW})}
Let $I=(f_{0},f_{1},...,f_{k}) \subseteq R$ be an equigenerated ideal of degree $r$. Then there exists a $D$ such that for all $d \geq D$, we have 
$$\beta_{i,j+rd}(I^d) \neq 0 \Longleftrightarrow \beta_{i,j+rD}(I^D) \neq 0.$$ 
\end{thm}

\subsection{Stabilization Index}
The power of $I$, $D$, for which we see a stabilized Betti table shape in Theorem\autoref{T: Theorem 2.9} is referred to as the \textit{stabilization index} of $I$, also defined by Whieldon.
\begin{defn}
[Definition 5.1 in \cite{GW}]
Let $I=(f_{0},f_{1},...,f_{k}) \subseteq R$ be an equigenerated ideal of degree $r$. Let the \textit{stabilization index Stab$(I)$} of $I$ be the smallest $D$ such that for all $d \geq D$,
$$\beta_{i,j+rd}(I^d) \neq 0 \Longleftrightarrow \beta_{i,j+rD}(I^D) \neq 0.$$ 
\end{defn}
Below is an example to aid in the understanding of these concepts.
\begin{example} \label{E: Example 2.11}
Let $I=(x_{1}x_{2}x_{3}x_{4},x_{2}^4,x_{1}x_{4}^3) \subseteq \mathbbm{k}[x_{1},x_{2},x_{3},x_{4}]$. If we consider the Betti tables of the first few powers of $I$,
$$
\begin{tabular}{cccccccc}
$I^1$ & & & $I^2$ & & & $I^3$ & \\
 &
\begin{tabular}{cccc}
- & 0 & 1 & 2\\
total: & 3 & 3 & 1\\
4: & 3 & $\cdot$ & $\cdot$\\
5: & $\cdot$ & 1 & $\cdot$\\
6: & $\cdot$ & 1 & $\cdot$\\
7: & $\cdot$ & 1 & 1\\
\end{tabular}
& & &
\begin{tabular}{cccc}
- & 0 & 1 & 2\\
total: & 6 & 9 & 4\\
8: & 6 & $\cdot$ & $\cdot$\\
9: & $\cdot$ & 3 & $\cdot$\\
10: & $\cdot$ & 4 & 2\\
11: & $\cdot$ & 2 & 2\\
\end{tabular}
& & &
\begin{tabular}{cccc}
- & 0 & 1 & 2\\
total: & 10 & 18 & 9\\
12: & 10 & $\cdot$ & $\cdot$\\
13: & $\cdot$ & 6 & $\cdot$\\
14: & $\cdot$ & 9 & 6\\
15: & $\cdot$ & 3 & 3\\
\end{tabular}\\
\end{tabular}
$$
we see two distinct Betti table shapes, first expressed in the Betti tables of $I^1$ and $I^2$. One can test higher powers of $I$ and see that it appears that $\beta_{i,j+4(2)}(I^2) \neq 0 \Longleftrightarrow \beta_{i,j+4(d)}(I^d) \neq 0$, for all $i$ and $j$, and $d \geq 2$. Given Theorem\autoref{T: Theorem 2.9}, it appears that all Betti tables of $I^d$, where $d \geq 2$, will have the same shape as $I^2$. Therefore, we predict that Stab$(I)=2$.
\end{example}
Though Whieldon's proof of Theorem\autoref{T: Theorem 2.9} ensures that the Stab$(I)$ is finite, there is no known general formula for the Stab$(I)$. Whieldon provides a conjecture for the formula for the stabilization index of edge ideals; see \cite{GW} Conjecture 5.2. Below we extend Whieldon's concept of the Stab$(I)$ to distinguish the different powers of $I$, $d$, where we see changes in the shape of the Betti table of $I^d$, up to the Stab$(I)$.
\subsection{Stabilization Sequence}
\begin{defn}
Let $I \subseteq R$. Let the \textit{stabilization sequence StabSeq$(I)$} of $I$ be the sequence of powers for which we see new shapes of the Betti tables of $I$.
$$\text{StabSeq}(I)=\{d : I^d\ \text{does not share the same Betti table shape as}\ I^{d-1}, d \in \mathbbm{Z}^+\}.$$
\end{defn}
In Example\autoref{E: Example 2.11}, we saw two distinct Betti table shapes, first expressed at $I^1$ and $I^2$, with a predicted Stab$(I)=2$. Therefore, $1,2 \in \text{StabSeq}(I)$ and $3 \not\in \text{StabSeq}(I)$. Our predicted stabilization sequence is therefore
$$
\text{StabSeq}(I)=\{1,2\}.
$$

By definition all stabilization sequences will contain 1, since 1 is the first power at which we can observe a Betti table of $I^d$. All stabilization sequences in this paper were determined using the Macaulay2 \cite{M2} package \texttt{StabSeq.m2}, created by the author. The source code for this package is displayed in \hyperref[Appendix B]{Section~\ref*{Appendix B}}.

\section{Questions and Results} \label{Q: Questions and Results}
There are many questions we can ask about the stabilization index and sequence of a given ideal. The rest of this paper focuses on exploring these questions through an example-based investigation.

\subsection{Question 1. \textnormal{\textit{If $\beta (I^d)$ and $\beta (I^{d+1})$ have the same shape, is Stab$(I)=d$?}}}\

One might hope that the shape of the Betti tables of $\{I^d\}$ will stabilize when consecutive powers share the same shape. That is, Stab$(I)=d$ where $d$ is the smallest power such that, for all $i$ and $j$, and some $r$,
$$
\beta_{i,j+rd}(I^d) \neq 0 \Longleftrightarrow \beta_{i,j+r(d+1)}(I^{d+1}) \neq 0.
$$

This would make determining the stabilization index very straightforward, for we would only have to check the Betti table shapes of powers of $I$ until consecutive powers share the same shape. This appears to be the case for the ideals in Example\autoref{E:Example 2.7} and Example\autoref{E: Example 2.11}. It also follows from Theorem 2.1 of \cite{EG} that the Betti tables of all complete intersection monomial equigenerated ideals of the form
$$
I=(x_{1}^s,x_{2}^s,\dots,x_{n}^s) \subseteq R
$$
will stabilize at $I^1$ and have the following shape, for all $d \in \mathbbm{Z^+}:$

$$
\begin{tabular}{cc}
$I^d$ & \\
&
\begin{tabular}{ccccccc}
- & 0 & 1 & $\cdots$ & $n-2$ & $n-1$\\
total: & $\ast$ & $\ast$ & $\cdots$ & $\ast$& $\ast$\\ \cline{2-2}
$ds$: & \multicolumn{1}{|c|}{$\ast$} & $\cdot$ & $\cdots$ & $\cdot$ & $\cdot$\\ \cline{2-3}
$(d+1)s-1$: & $\cdot$ & \multicolumn{1}{|c|}{$\ast$} & $\cdots$ & $\cdot$ & $\cdot$\\ \cline{3-3}
$\vdots$ & $\cdot$ & $\cdot$ & $\ddots$ & $\cdot$ & $\cdot$\\ \cline{5-5}
$(d+(n-2))s-(n-2)$: & $\cdot$ & $\cdot$ & $\cdots$ & \multicolumn{1}{|c|}{$\ast$} & $\cdot$\\ \cline{5-6}
$(d+(n-1))s-(n-1)$: & $\cdot$ & $\cdot$ & $\cdots$ & $\cdot$ & \multicolumn{1}{|c|}{$\ast$}\\ \cline{6-6}
\end{tabular}\\
\end{tabular}
$$

Therefore the ideal in Example\autoref{E:Example 2.7}, $I=(x_{1}^2,x_{2}^2,x_{3}^2,x_{4}^2) \subseteq k[x_{1},x_{2},x_{3},x_{4}]$, will have the following Betti table shape, for all $d \in \mathbbm{Z^+}$,
$$
\begin{tabular}{cc}
$I^d$ & \\
&
\begin{tabular}{cccccc}
- & 0 & 1 & 2 & 3\\
total: & $\ast$ & $\ast$ & $\ast$ & $\ast$\\ \cline{2-2}
$2d$: & \multicolumn{1}{|c|}{$\ast$} & $\cdot$ & $\cdot$ & $\cdot$\\ \cline{2-3}
$2d+1$: & $\cdot$ & \multicolumn{1}{|c|}{$\ast$} & $\cdot$ & $\cdot$\\ \cline{3-4}
$2d+2$: & $\cdot$ & $\cdot$ & \multicolumn{1}{|c|}{$\ast$} & $\cdot$\\ \cline{4-5}
$2d+3$: & $\cdot$ & $\cdot$ & $\cdot$ & \multicolumn{1}{|c|}{$\ast$}\\ \cline{5-5}
\end{tabular}\\
\end{tabular}
$$
and stabilization index and sequence
$$
\begin{tabular}{ccc}
Stab$(I)=1$, & StabSeq$(I)=\{1\}.$\\
\end{tabular}
$$

However, the following example shows that we are not guaranteed Betti table shape stabilization when consecutive powers share the same shape.
\begin{example}  \label{E:Example 3.1}
Let $I=(x_{1}^3x_{2},x_{2}^4,x_{1}^2x_{3}^2,x_{2}^3x_{3}) \subseteq \mathbbm{k}[x_{1},x_{2},x_{3}]$. Below are the Betti tables of the first few powers of $I$.
$$
\begin{tabular}{cccccccc}
$I^1$ & & & $I^2$ & & & $I^3$ & \\
 &
\begin{tabular}{cccc}
- & 0 & 1 & 2 \\
total: & 4 & 5 & 2 \\
4: & 4 & 1 & $\cdot$ \\
5: & $\cdot$ & 1 & $\cdot$ \\
6: & $\cdot$ & 3 & 2 \\
\end{tabular}
& & &
\begin{tabular}{cccc}
- & 0 & 1 & 2 \\
total: & 10 & 15 & 6 \\
8: & 10 & 5 & $\cdot$ \\
9: & $\cdot$ & 3 & $\cdot$ \\
10: & $\cdot$ & 7 & 6 \\
\end{tabular}
& & &
\begin{tabular}{cccc}
- & 0 & 1 & 2 \\
total: & 20 & 32 & 13 \\
12: & 20 & 14 & 1 \\
13: & $\cdot$ & 7 & 1  \\
14: & $\cdot$ & 11 & 11 \\
\end{tabular}\\
\end{tabular}
$$

In the above tables we see that $I^1$ and $I^2$ share the same shape, that is, for all $i$ and $j$,
$$
\beta_{i,j+4(1)}(I^1) \neq 0 \Longleftrightarrow \beta_{i,j+4(2)}(I^2) \neq 0.
$$
However, the Stab$(I) \neq 1$ since the Betti table of $I^3$ differs in shape than that of $I^1$. When we check the shape of the Betti tables of $I^d$ where $d \geq 3$ we see that, for all $i$ and $j$,
$$
\beta_{i,j+4(3)}(I^3) \neq 0 \Longleftrightarrow \beta_{i,j+4(d)}(I^d) \neq 0.
$$
Therefore, we predict the stabilized shape of the Betti tables of $I^d$ for $d \geq 3$ to be

$$
\begin{tabular}{cc}
$I^d$ & \\
&
\begin{tabular}{cccc}
- & 0 & 1 & 2 \\
total: & $\ast$ & $\ast$ & $\ast$ \\ \cline{2-4}
$4d$: & \multicolumn{1}{|c}{$\ast$} & $\ast$ &  \multicolumn{1}{c|}{$\ast$} \\ \cline{2-2}
$4d+1$: & $\cdot$ & \multicolumn{1}{|c}{$\ast$} & \multicolumn{1}{c|}{$\ast$}  \\
$4d+2$: & $\cdot$ &  \multicolumn{1}{|c}{$\ast$} & \multicolumn{1}{c|}{$\ast$} \\ \cline{3-4}
\end{tabular}\\
\end{tabular}
$$
with the following stabilization index and sequence 
$$
\begin{tabular}{ccc}
Stab$(I)=3$ & & StabSeq$(I)=\{1,3\}.$\\
\end{tabular}
$$
\end{example}

From this example, we know that we are not guaranteed to have found the stabilized Betti table shape, of a given ideal, when the Betti tables of two consecutive powers share the same shape. In the many cases where we do not have the tools to generalize the Betti tables of a given ideal, we have to make hypothesis for the stabilization index. 

\subsection{Question 2. \textnormal{\textit{Does there exist some $x$ such that if $\beta_{i,j+r(D)}(I^D) \neq 0 \Longleftrightarrow \beta_{i,j+r(d)}(I^{d}) \neq 0,$ for all $d \in \{(D+1),\dots,(D+x)\},$ then Stab$(I) \leq D$?}}}\

Can we be confident that the shape of the Betti tables of $\{I^d\}$ stabilize when a large enough number of consecutive Betti tables in this collection share the same shape? Determining a bound on this $x$ would allow for far more straightforward computational investigation into stabilization indices and sequences, as there would be a recognizable point which would mark the stabilization. It is likely that this $x$ is dependant upon characteristics of the ideal in question, however what specific characteristics these may be we cannot answer here. From the following example we know that such an $x$ must be at least 7.

\begin{example} \label{E: Example 3.4}
Let $I=(x_{1}^3x_{2}^3x_{3}^3,x_{2}^6x_{3}^3,x_{1}^4x_{3}^5,x_{1}^8x_{2}) \subseteq \mathbbm{k}[x_{1},x_{2},x_{3}]$. Below are the Betti tables of $I^d$ for which $I^d$ does not share the same Betti table shape as $I^{d-1}$.

$$
\begin{tabular}{cccccccc}
$I^1$ & & & $I^2$ & & & $I^3$ &\\
 &
\begin{tabular}{cccc}
- & 0 & 1 & 2 \\
total: & 4 & 4 & 1 \\
9: & 4 & $\cdot$ & $\cdot$ \\
10: & $\cdot$ & $\cdot$ & $\cdot$ \\
11: & $\cdot$ & 2 & $\cdot$ \\
12: & $\cdot$ & $\cdot$ & $\cdot$ \\
13: & $\cdot$ & 2 & $\cdot$ \\
14: & $\cdot$ & $\cdot$ & 1 \\
\end{tabular}
& & &
\begin{tabular}{cccc}
- & 0 & 1 & 2 \\
total: & 10 & 15 & 6 \\
18: & 10 & $\cdot$ & $\cdot$ \\
19: & $\cdot$ & 1 & $\cdot$ \\
20: & $\cdot$ & 8 & $\cdot$ \\
21: & $\cdot$ & $\cdot$ & 1 \\
22: & $\cdot$ & 6 & 2 \\
23: & $\cdot$ & $\cdot$ & 3 \\
\end{tabular}
& & &
\begin{tabular}{cccc}
- & 0 & 1 & 2 \\
total: & 20 & 35 & 16 \\
27: & 20 & $\cdot$ & $\cdot$ \\
28: & $\cdot$ & 4 & $\cdot$ \\
29: & $\cdot$ & 21 & 1 \\
30: & $\cdot$ & $\cdot$ & 5 \\
31: & $\cdot$ & 10 & 5 \\
32: & $\cdot$ & $\cdot$ & 5 \\
\end{tabular}\\
\end{tabular}
$$

$$
\begin{tabular}{cccccccc}
$I^4$ & & & $I^6$ & & & $I^{13}$ &\\
 &
\begin{tabular}{cccc}
- & 0 & 1 & 2 \\
total: & 35 & 66 & 32 \\
36: & 35 & 1 & $\cdot$ \\
37: & $\cdot$ & 10 & $\cdot$ \\
38: & $\cdot$ & 41 & 4 \\
39: & $\cdot$ & $\cdot$ & 13 \\
40: & $\cdot$ & 14 & 8 \\
41: & $\cdot$ & $\cdot$ & 7 \\
\end{tabular}
& & &
\begin{tabular}{cccc}
- & 0 & 1 & 2 \\
total: & 84 & 166 & 83 \\
54: & 84 & 10 & $\cdot$ \\
55: & $\cdot$ & 39 & 1 \\
56: & $\cdot$ & 95 & 20 \\
57: & $\cdot$ & $\cdot$ & 37 \\
58: & $\cdot$ & 22 & 14 \\
59: & $\cdot$ & $\cdot$ & 11 \\
\end{tabular}
& & &
\begin{tabular}{cccc}
- & 0 & 1 & 2 \\
total: & 560 & 1101 & 542 \\
117: & 560 & 255 & 1 \\
118: & $\cdot$ & 430 & 130 \\
119: & $\cdot$ & 336 & 216 \\
120: & $\cdot$ & $\cdot$ & 135 \\
121: & $\cdot$ & 50 & 35 \\
122: & $\cdot$ & $\cdot$ & 25 \\
\end{tabular}\\
\end{tabular}
$$

We have checked the Betti tables shapes of $I^d$ for $d \leq 100$ and found the same Betti table shape as that of $I^{13}$. Thus we hypothesize the stabilization index and sequence of $I$ to be
$$
\begin{tabular}{ccc}
Stab$(I)=13$ & & StabSeq$(I)=\{1,2,3,4,6,13\}.$\\
\end{tabular}
$$
\end{example}
From this example, which is derived from a relatively simple homogeneous ideal, we see a more complicated stabilization sequence, one which contains 6 elements, and a relatively large stabilization index occurring at $I^d$ where $d>10$. Another interesting thing to note from this example is that the Betti table of $I^d$ appeared to have stabilized at $I^6$, that is until we checked $I^{6+x}$ where $x=7$. Though checking the Betti table shape 7 powers higher after a new Betti table shape is not computationally intensive, this example shows how little we know about the bound on $x$. For if such a simple ideal can produce a relatively complicated stabilization index and sequence, with a bound on $x\geq 7$, what can we expect from more complicated ideals? We expect that ideals with a greater number of and degree of variables and generators may produce these more complicated stabilization indices and sequences. 

\subsection{Question 3. \textnormal{\textit{In what ways can we characterize stabilization sequences?}}} \label{Q: Question 3}\

We acknowledge that this question is quite general and could be used to guide investigations of families of ideals $\{I_{i}\}$ whose Betti table shape does not necessarily stabilize. We focus on exploring this question for collections of ideals $\{I^d\}$ where $I$ is an equigenerated ideal, of fixed degree or a fixed number of generators. In Example\autoref{E: Example 3.4}, we saw an ideal with a stabilization sequence containing 6 elements, whereas all the previous examples contained at most 2 elements. Given these examples we may wonder about the cardinality of other stabilization sequences. Are there ideals whose stabilization sequences are of arbitrarily large cardinality? We may also ask whether some stabilization sequence are polynomial? The following example presents a surprising stabilization sequence.

\begin{example} \label{E: Example 3.6}
Let $I=(x_{1}^8x_{2}^8x_{3}^8,x_{1}^4x_{2}^{16}x_{3}^4,x_{1}x_{2}^{23},x_{2}^{12}x_{3}^{12}) \subseteq \mathbbm{k}[x_{1},x_{2},x_{3}]$. One can test shapes of the Betti tables of $I^d$ where $1 \leq d \leq 20$ and find the following stabilization sequence,
$$
\text{StabSeq}(I)=\{1,2,3,5,7,9,15,17,19\}.
$$
\end{example}

We have checked the Betti tables shapes of $I^d$ for $d \leq 100$ and found the same Betti table shape as that of $I^{19}$. Though we do not see a stabilization sequence of arbitrarily large cardinality for this ideal, we do see a stabilization sequence containing far more elements than previous examples, $|\text{StabSeq}(I)|=9$. In this example we also see rather disconcerting behaviour of the elements in the stabilization sequence of $I$, for elements $3,5,7,9$ grow by increments of 2 yet the following element is 15. This behaviour is problematic because it leads us to assume that elements will continue to appear by increments of 2, then produces an element incremented by 6. This example provides further evidence that we should be cautious when dealing with stabilization sequences of ideals. For in many cases where we do not currently have the tools to definitively state them, there is no guarantee that any structure or pattern they appear to exhibit will continue. 

\subsection{Question 4. \textnormal{\textit{For equigenerated ideals, of fixed degree or a fixed number of generators, can the stabilization index be arbitrarily large?}}} \label{Q: Question 4}\

In \autoref{Q: Question 3} we ask whether there exists stabilization sequence of arbitrarily large cardinality. A natural extension of this question is whether there exists equigenerated ideals, of fixed degree or a fixed number of generators, whose stabilization index is arbitrarily large. In Example\autoref{E: Example 3.4} and Example\autoref{E: Example 3.6}, we saw stabilization indices of 11 and 19 respectively. One might expect, given the relative simplicity of the ideal in these examples, that more complicated ideals, will stabilize at higher powers, perhaps even some which stabilize at $I^d$ where $d>>0$. We predict that the existence of such indices depend upon the number of and degree of the variables and generators of the ideal. 

\subsection{Linearly connected family of ideals}
Until now we have been working with collections of ideals $\{I^d\}$ as we vary $d$. These collections are examples of \textit{graded families} of ideals. 

\begin{defn}
A collection of ideals  $\{I_{i}\} \subseteq R$ is a \textit{graded family} if, for all $i$ and $j$, $I_{i}I_{j} \subset I_{i+j}$.
\end{defn}

It is understood that there is inherent structure between ideals in a graded family. We now define a \textit{linearly connected family} of ideals, which are not a graded family in general. Recall we may write a monomial in $R$ as $x^{\alpha}=x_{1}^{{\alpha}_{1}}x_{2}^{{\alpha}_{2}}\cdots x_{n}^{{\alpha}_{n}}$.

\begin{defn}
A collection of monomial ideals $\{I_{j}\} \subseteq R$ is a \textit{linearly connected family} if there exists linear functions $\alpha_{k}^i(j)$ such that
$$
I_{j}=(x^{\alpha^1(j)},x^{\alpha^2(j)},...,x^{\alpha^l(j) })\ \text{where}\ \alpha^i(j)=(\alpha_{1}^i(j),\alpha_{2}^i(j),...,\alpha_{n}^i(j))
$$
$\text{for}\ 1 \leq i \leq l\ \text{and}\  j \in \mathbbm{Z}^+.$
\end{defn}

\begin{example}
The homogeneous ideal $\{I_{b}\}=(x_{1}^bx_{2}^{3b},x_{2}^{4b},x_{1}^{2b}x_{3}^{2b}) \subseteq \mathbbm{k}[x_{1},x_{2},x_{3}]$ is a \textit{linearly connected family} as the powers of the variables of the generators of $I$ are all linear functions dependant on $b$. The $\alpha^i(b)$ are listed below
$$
\begin{tabular}{ccc}
$\alpha^1(b)=(b,3b,0)$, & $\alpha^2(b)=(0,4b,0)$, & $\alpha^3(b)=(2b,0,2b)$.\\
\end{tabular}
$$
\end{example}

We will now discuss the linearly connected family of homogeneous ideals $\{I_{n}\}$, where 
$$I_{n}=(a^{2n}b^{2n}c^{2n},b^{4n}c^{2n},a^{3n}c^{3n},a^{6n-1}b)\subseteq \mathbbm{k}[a,b,c],$$
as we vary $n$ and the powers of $I_{n}$. For $n=1$, $I_{1}=(a^2b^2c^2,b^4c^2,a^3c^3,a^5b)$, we observe three distinct Betti table shapes. Below are the Betti tables of $I_{1}^d$ for which $I_{1}^d$ does not share the same Betti table shape as $I_{1}^{d-1}$.

$$
\begin{tabular}{cccccccc}
$I_{1}^1$ & & & $I_{1}^2$ & & & $I_{1}^6$\\
 & 
 \begin{tabular}{cccc}
- & 1 & 2 & 3 \\
total: & 4 & 4 & 1 \\
4: & 4 & $\cdot$  & $\cdot$ \\
5: & $\cdot$ & 2 & $\cdot$ \\
6: & $\cdot$ & 2 & 1 \\
\end{tabular}
& & &
\begin{tabular}{cccc}
- & 1 & 2 & 3 \\
total: & 10 & 15 & 1 \\
11: & 10 & 1  & $\cdot$ \\
12: & $\cdot$ & 8 & 1 \\
13: & $\cdot$ & 6 & 5 \\
\end{tabular}
 & & &
\begin{tabular}{cccc}
- & 1 & 2 & 3 \\
total: & 84 & 162 & 79 \\
35: & 84 & 49  & 1 \\
36: & $\cdot$ & 91 & 53 \\
37: & $\cdot$ & 22 & 25 \\
\end{tabular}\\
\end{tabular}
$$

The above result is not surprising, especially when considering the stabilization index and sequences of the ideal in Example\autoref{E: Example 3.4} and Example\autoref{E: Example 3.6}. However, when we increase the value of $n$ and test the resulting stabilization sequences of $I_{n}$, we see very surprising results. Below we have listed the stabilization sequences of $I_{n}$ for $ 1 \leq n \leq 8$, determined by testing the shapes of the Betti tables of $I_{n}^d$ for $d\leq100$.
$$
\begin{tabular}{l}
$\text{StabSeq}(I_{1})= \{1,2,6\}$\\
 \\
$\text{StabSeq}(I_{2})= \{1,2,3,5,6,11\}$\\
 \\
$\text{StabSeq}(I_{3})=\{1,2,3,5,6,11,17,23\}$\\
 \\
$\text{StabSeq}(I_{4})=\{1,2,3,5,6,11,17,23,29,35\}$\\
 \\
$\text{StabSeq}(I_{5})=\{1,2,3,5,6,11,17,23,29,35,41,47\}$\\
 \\
$\text{StabSeq}(I_{6})=\{1,2,3,5,6,11,17,23,29,35,41,47,53,59\}$\\
 \\
$\text{StabSeq}(I_{7})=\{1,2,3,5,6,11,17,23,29,35,41,47,53,59,65,71\}$\\
 \\
$\text{StabSeq}(I_{8})=\{1,2,3,5,6,11,17,23,29,35,41,47,53,59,65,71,77,83\}$\\
\end{tabular}
$$

We had no reason to expect that the stabilization indices or sequences of a linearly connected family of ideals should have any relation. Yet we see a very clear pattern and apparent structure arising in this particular linearly connected family of homogeneous ideals. We see that the stabilization index of $I_{n}$, for $n \geq 2$, appears to be given by the linear function
$$
\text{Stab}(I_{n})=12n-13
$$
and the stabilization sequence of $I_{n}$ appears to be given by
$$
\text{StabSeq}(I_{n})=\text{StabSeq}(I_{n-1}) \cup \{12n-13,12n-19\}.
$$

These are entirely unexpected relationships, which raise more questions than answers. In \autoref{Q: Question 3} we ask in what ways can we characterize the stabilization sequences? From this related collection we see that there are in fact stabilization sequences with potentially large cardinality. For the cardinality of $I_{n}$ for $n\geq2$ appears to be given by the linear function
$$
|\text{StabSeq}(I_{n})|=2n+2.
$$

We also see that the stabilization sequences appear to be somewhat linear, that is, the stabilization sequence of $I_{n}$ not only appears to be a union of the stabilization sequence of $I_{n-1}$ and $\{12n-13,12n-19\}$ but also expressible as 
$$
\text{StabSeq}(I_{n})=\{1,2,3,5,6\} \cup \{d:d=\{11+6m\}\ \text{for all}\ 0\leq m\leq n-1, m \in \mathbbm{Z^+}\}
$$
for $n\geq3$. In \autoref{Q: Question 4} we ask whether there exists a collection of equigenerated ideals, of fixed degree or fixed number of generators, whose Betti tables stabilized in shape at an arbitrarily large power? If the structure in this related collection continues, it appears that as $n$ becomes arbitrarily large, so will the stabilization index of $I_{n}$, being dependant upon
$$
\text{Stab}(I_{n})=12n-13.
$$

\section{Future Research}
The unexpected structure of the stabilization indices and sequences found in the linearly related family of homogeneous ideals above presents a wealth of questions for future research. Below we list some questions which we believe to be fruitful avenues for future research.
\begin{question}
Do all linearly related family of homogeneous ideals exhibit structure in their stabilization indices and sequences? If so, in what ways can we characterize this structure?
\end{question}
\begin{question}
What sorts of functions are the stabilization indices and sequences of other linearly related family of homogeneous ideals dependant upon? Are these functions always linear?
\end{question}
\begin{question}
Is there a formula for Stab$(I)$ and StabSeq$(I)$ for certain classes of ideals? If so what characteristics of the ideals does it depend upon?
\end{question}

$$
$$
$$
$$
$$
$$
$$
$$

\break
\section{Appendix A.} \label{Appendix A}

Macaulay2 \cite{M2} commands used in this paper.
$$
\begin{tabular}{ll} \cline{1-2}
\multicolumn{1}{|c|}{Operation} & \multicolumn{1}{|c|}{M2 Command}\\ \cline{1-2}
\multicolumn{1}{|l|}{Defining Field:} & \multicolumn{1}{|l|}{
\begin{minipage}{2.9in}
\begin{verbatim}
R=QQ[input variables]
\end{verbatim}
\end{minipage}}
\\
\multicolumn{1}{|l|}{Defining ideal:} & \multicolumn{1}{|l|}{
\begin{minipage}{2.9in}
\begin{verbatim}
I=ideal(input generators of ideal)
\end{verbatim}
\end{minipage}}
\\
\multicolumn{1}{|l|}{Define ideal as a module:} & \multicolumn{1}{|l|}{
\begin{minipage}{2.9in}
\begin{verbatim}
M=module I
\end{verbatim}
\end{minipage}}
\\
\multicolumn{1}{|l|}{Define resolution of ideal:} & \multicolumn{1}{|l|}{
\begin{minipage}{2.9in}
\begin{verbatim}
C=res minimalPresentation M
\end{verbatim}
\end{minipage}}
\\ 
\multicolumn{1}{|l|}{View resolution of ideal:} & \multicolumn{1}{|l|}{
\begin{minipage}{2.9in}
\begin{verbatim}
C.dd
\end{verbatim}
\end{minipage}}
\\ 
\multicolumn{1}{|l|}{View Betti table of ideal:} & \multicolumn{1}{|l|}{
\begin{minipage}{2.9in}
\begin{verbatim}
betti res minimalPresentation M
\end{verbatim}
\end{minipage}}
\\
\multicolumn{1}{|l|}{Load stabilization sequence package from} & \multicolumn{1}{|l|}{
\texttt{load "StabSeq.m2"}}
\\
\multicolumn{1}{|l|}{Macaulay2 code/ directory:} & \multicolumn{1}{|l|}{}
\\
\multicolumn{1}{|l|}{Determine stabilization sequence of ideal:} & \multicolumn{1}{|l|}{
\texttt{StabSeq(I, input max power to check,}}\\
\multicolumn{1}{|l|}{} & \multicolumn{1}{|l|}{
\texttt{IncludeBettis => true)}}
\\  \cline{1-2}
\end{tabular}
$$

Below we provide the Macaulay2 \cite{M2} session used in finding the minimal graded free resolution and Betti table of the ideal in Example\autoref{E: Example 2.1} and Example\autoref{E:Example 2.3} with commentary. We begin by defining \texttt{R} as the polynomial ring that we will be working in. We then construct the ideal \texttt{I}, in \texttt{R}, and define \texttt{M} as the module of \texttt{I}.\\
\begin{verbatim}
i1 : R=QQ[x_1,x_2,x_3]

o1 = R

o1 : PolynomialRing

i2 : I=ideal(x_1*x_2^2,x_1*x_3^2,x_2^3,x_1^3)

               2     2   3   3
o2 = ideal (x x , x x , x , x )
             1 2   1 3   2   1

o2 : Ideal of R

i3 : M=module I

o3 = image | x_1x_2^2 x_1x_3^2 x_2^3 x_1^3 |

                             1
o3 : R-module, submodule of R
\end{verbatim}

$ $

Given the module of \texttt{I}, \texttt{M}, we define \texttt{C} as the minimally graded free resolution of \texttt{M}. The command \texttt{C.dd} allows us to view the minimally graded free resolution of \texttt{M} in its entirety.\\
\break
\begin{verbatim}
i4 : C=res minimalPresentation M

      4      4      1
o4 = R  <-- R  <-- R  <-- 0
                           
     0      1      2      3

o4 : ChainComplex

i5 : C.dd

          4                                         4
o5 = 0 : R  <------------------------------------- R  : 1
               {3} | -x_2 x_1^2  0      -x_3^2 |
               {3} | 0    0      x_1^2  x_2^2  |
               {3} | x_1  0      0      0      |
               {3} | 0    -x_2^2 -x_3^2 0      |

          4                      1
     1 : R  <------------------ R  : 2
               {4} | 0      |
               {5} | x_3^2  |
               {5} | -x_2^2 |
               {5} | x_1^2  |

          1
     2 : R  <----- 0 : 3
               0

o5 : ChainComplexMap
\end{verbatim}

$ $

Lastly, we view the Betti table of the module of \texttt{I}, \texttt{M}.\\
\begin{verbatim}
i6 : betti res minimalPresentation M

            0 1 2
o6 = total: 4 4 1
         3: 4 1 .
         4: . 3 .
         5: . . 1

o6 : BettiTally
\end{verbatim}

$ $

Below we provide the Macaulay2 \cite{M2} session used in finding the stabilization sequence of the ideal in Example\autoref{E:Example 3.1} and Example\autoref{E: Example 3.6} with commentary. Similar to the previous session, we begin by defining \texttt{R} as the polynomial ring then construct the ideal \texttt{I}, in \texttt{R}.\\
\begin{verbatim}
i1 : R=QQ[x_1,x_2,x_3]

o1 = R

o1 : PolynomialRing

i2 : I=ideal(x_1^3*x_2,x_2^4,x_1^2*x_3^2,x_2^3*x_3)

             3     4   2 2   3
o2 = ideal (x x , x , x x , x x )
             1 2   2   1 3   2 3

o2 : Ideal of R
\end{verbatim}

$ $

We load the stabilization sequence package, provided in \hyperref[Appendix B]{Section~\ref*{Appendix B}}, from the /Library/Application Support/Macaulay2/code/ directory.\\
\begin{verbatim}
i3 : load "StabSeq.m2"
\end{verbatim}

$ $

For full documentation of the stabilization package, see \hyperref[Appendix B]{Section~\ref*{Appendix B}}. The stabilization sequence package computes the stabilization sequence of the inputed ideal (\texttt{I} in the following example) up to the indicated maximum power (\texttt{25} in the following example). By including the optional argument \texttt{IncludeBettis => true} the \texttt{StabSeq} function also outputs the Betti tables for which we see a new Betti table shape. Below is the stabilization sequence and Betti tables found in Example\autoref{E:Example 3.1}.\\
\begin{verbatim}
i4 : StabSeq(I,25,IncludeBettis => true)
        0 1 2          0  1  2
{total: 4 5 2, total: 20 32 13}
     4: 4 1 .     12: 20 14  1
     5: . 1 .     13:  .  7  1
     6: . 3 2     14:  . 11 11
o4 = {1, 3}

o4 : List
\end{verbatim}

$ $

Without the optional argument \texttt{IncludeBettis => true} the \texttt{StabSeq} function simply outputs the stabilization sequence of the inputed ideal (\texttt{J} in the following example) up to the indicated maximum power (\texttt{25} in the following example). Below is the stabilization sequence and Betti tables found in Example\autoref{E: Example 3.6}.\\
\begin{verbatim}
i5 : J=ideal(x_1^8*x_2^8*x_3^8,x_1^4*x_2^(16)*x_3^4,x_1*x_2^(23),x_2^(12)*x_3^(12))

             8 8 8   4 16 4     23   12 12
o5 = ideal (x x x , x x  x , x x  , x  x  )
             1 2 3   1 2  3   1 2    2  3

o5 : Ideal of R

i6 : StabSeq(J,25)

o6 = {1, 2, 3, 5, 7, 9, 15, 17, 19}

o6 : List
\end{verbatim}

\section{Appendix B.} \label{Appendix B}

Below we provide the source code for \texttt{StabSeq.m2}. To use this algorithm, save it as a .m2 file in /Library/Application Support/Macaulay2/code/. Note that the \texttt{--} in the source code below corresponded to comments to aid in the understanding of this package.\\

\begin{lstlisting}
StabSeq =  method(Options => {IncludeBettis => false})
needsPackage "BoijSoederberg"

-- By Aaron Slobodin, September 2017, open source.

-- Algorithm returns the Stabilization Sequence (up to inputed max power) of the inputed ideal.
-- The algorithm includes the optional argument to print the Betti Tables of all powers of the inputed ideal that are found in the determined Stabilization Sequence.
-- Note that the Stabilization Sequence of a given ideal, I, is {i+1: I^(i+1) does not share the same Betti table shape as I^i, i in ZZ+}

StabSeq(Ideal,ZZ) := o -> (I, PowerCap) -> (
	StabSeqList := {1};
	CurrentBetti := betti res minimalPresentation module I;
	BettiLister := {CurrentBetti};
	CurrentMatrix = matrix CurrentBetti;
	for i from 1 to (PowerCap-1) do (
		NewBetti := betti res minimalPresentation module I^(i+1);
		NewMatrix := matrix NewBetti;
		
-- Check to see if dimensions of the Betti tables of I^i and I^(i+1) are the same.

		if (numgens source CurrentMatrix - 1) != (numgens source NewMatrix - 1) or (numgens target CurrentMatrix - 1) != (numgens target NewMatrix - 1) then (
		
-- If they dimensions are not equal, they have different Betti table shapes, resulting in
-- 1. the stabilization sequence (StabSeqList) gaining the element i+1.
-- 2. the list of Betti tables (BettiLister) gains the element NewBetti

				StabSeqList = append(StabSeqList, i+1);
				BettiLister = append(BettiLister, NewBetti);
		)
		
-- Check to see if Betti table of I^i and I^(i+1) exhibit the same shape.

		else (
			for j from 0 to (numgens source CurrentMatrix - 1) do (
				for k from 0 to (numgens target CurrentMatrix - 1) do (
					if CurrentMatrix_j_k != 0 and NewMatrix_j_k == 0 or CurrentMatrix_j_k == 0 and NewMatrix_j_k != 0 then (
					
-- If the Betti tables of I^i and I^i+1 differ in shape then
-- 1. the stabilization sequence (StabSeqList) gaining the element i+1.
-- 2. the list of Betti tables (BettiLister) gains the element NewBetti

						StabSeqList = append(StabSeqList, i+1);
						BettiLister = append(BettiLister, NewBetti);
					);
				);
			);
		);
		
-- Either the Betti table of I^i and I^(i+1) shared the same shape or differed in shape. Regardless, the new reference Betti table should therefore be the Betti table of I^(i+1).		

	CurrentMatrix = NewMatrix;
	);
	
-- Optional command to print the Betti tables of elements in the determined stabilization sequence.

if o.IncludeBettis then (<< unique BettiLister);

-- Return determined stabilization sequence.


return unique StabSeqList;
);
\end{lstlisting}

\section{Acknowledgments}

I would like to thank the Quest Summer Fellows Committee for providing me the funding and opportunity for this research and Dr. Sarah Mayes-Tang for her endless support as my host faculty advisor. Calculations in this paper were performed using the computer software Macaulay2 \cite{M2}.

\bibliography{SFBib}
\bibliographystyle{amsalpha}
\nocite{*}
\end{document}